\documentclass[12pt]{article}

\usepackage{amsmath}
\usepackage{amssymb}
\usepackage{amsthm}
\usepackage{graphicx}
\usepackage{fullpage}
\usepackage{textcomp}
\usepackage{cite}
\usepackage{url}
\usepackage{epstopdf}

\newtheorem{thm}{Theorem}
\newtheorem{lem}[thm]{Lemma}
\newtheorem{cor}[thm]{Corollary}

\theoremstyle{definition}
\newtheorem{defn}{Definition}
\newtheorem{rem}{Remark}

\newtheorem{as}{Assumption}

\title{\Large \bf On the Design of Campus Parking Systems with QoS guarantees\footnote{This work was supported in part by Science Foundation Ireland (SFI) grant 11/PI/1177 and in part by EU FP7 project INSIGHT under grant
    318225.}}

\author{Wynita Griggs\footnote{W. Griggs and R. Shorten are with the Hamilton Institute, National University of Ireland Maynooth, Maynooth, Co. Kildare, Ireland.  Corresponding author:  Wynita Griggs.  Phone:  +353-(0)1-7086100.  Fax:  +353-(0)1-7086269.  Email:  {\tt\small wynita.griggs@nuim.ie}}, Jia Yuan Yu\footnote{J. Y. Yu and R. Shorten are with IBM Research Ireland, Dublin, Ireland.}, Fabian Wirth\footnote{F. Wirth is with the University of Passau, Faculty of Computer Science and Mathematics, Innstrasse 33, 94032 Passau, Germany.}, Florian H\"{a}usler\footnote{F. H\"{a}usler is with Fraunhofer FOKUS, Berlin, Germany.}\\and Robert Shorten\addtocounter{footnote}{-3}\footnotemark[\value{footnote}] \addtocounter{footnote}{1}\footnotemark[\value{footnote}]}

\renewcommand{\P}{{\mathbb P}}
\newcommand{\R}{\mathbb{R}}
\newcommand{\N}{\mathbb{N}}
\usepackage{color}

\DeclareMathOperator*{\minimize}{minimize} 

\begin{document}

\maketitle
\pagestyle{plain}

\begin{abstract}
Parking spaces are resources that can be pooled together and shared, especially when there are complementary day-time and night-time users. We answer two design questions. First, given a quality of service requirement, how many spaces should be set aside  as contingency during day-time for night-time users?
Next, how can we replace the first-come-first-served access method by one that aims at optimal efficiency while keeping user preferences private?
\end{abstract}

\section{Introduction}\label{Introduction}

It was recently reported that over one year in a small Los Angeles business district, cars cruising for parking burned 47,000 gallons of gasoline and produced 730 tons of carbon dioxide \cite{shoup2007}.  Meanwhile, the consulting firm McKinsey recently claimed that the average car owner in Paris spends four years of his or her life searching for a parking space \cite{mckinsey}. The parking assignment problem associated with electric vehicles becomes even more acute. Due to the limited range of these vehicles, the marginal cost of expending energy to search for spaces may, in some cities, be prohibitively high. Thus, there is a real and compelling societal and economic need to revisit parking.

\bigskip

Increases in car ownership, inadequate public and private parking facilities, and synchronised demand, have led to serious temporal mismatches  in parking space supply and demand. Finding a parking space at certain times of the day is a non-trivial challenge. Furthermore, finding spaces is not only costly in terms of wasted time for the driver, but it also creates congestion and pollution affecting everyone. Thus, improving parking improves both economic efficiency and environmental quality.

\bigskip

% Recent advances in automotive technology, infrastructure deployment, and in the mathematics of large-scale optimisation and distributed congestion control, have opened up a plethora of new applications in the areas of {\em intelligent transportation systems} and {\em smarter cities applications}. Within this context, the issue of parking, which is both practically important and mathematically challenging to address, has been the focus of much of the recent interest from both practitioners and academics alike, giving rise to a number of new products from high profile companies.

% \bigskip

Indeed, already major companies and cities are responding to these challenges.  For example, many city authorities, in order to safeguard resident on-street parking from people commuting to a city through the day, have initiated schemes such as resident permit parking, and a number of commercial initiatives have emerged in the parking area.  SFpark and JustPark provide examples of city authorities and companies investing heavily in parking research and products within a smart cities context.

\bigskip

In this work, our objective is not to revisit prior work on the topic. Rather, we propose a new solution that is widely applicable in cities. 
\emph{We consider two nearby entities that have complementary supply and demand, i.e., there is a shortfall of parking space at one and a surplus at the other.} We call the entity with a shortfall of parking a {\em mini-city}. Many examples come to mind: university campuses, technology parks, shopping centers, and government facilities located in the suburbs.  In our work, we use the university campus as an example.

\bigskip

The second entity is the area surrounding the mini-city.
During business hours, while parking spaces are limited in the mini-city, many nearby residents leave home for work. Thus, there is an opportunity for the mini-city to use these vacated parking spaces, as well as from schools, hotels, apartment complexes and so on. Consequently, the mini-city has access to two classes of parking spaces:  \emph{premium} spaces, e.g., those located on the university campus; and \emph{secondary} spaces located nearby and perhaps serviced by a shuttle. Note that the secondary spaces are not owned by the mini-city but rather are rented from secondary parties or \emph{landlords}. Given this basic scheme, we consider two specific design issues. First, we wish to guarantee a quality of service for the landlords by setting aside reserve spaces in the mini-city as contingency for events where secondary spaces are suddenly unavailable. Second, we wish to ensure that the premium spaces are allocated optimally among \emph{users} (drivers) while preserving each user's privacy.

% \bigskip

% The objective in this paper is to show how such schemes can be designed in a manner that preserves an individual's privacy. We illustrate our designs with some simple examples.

%\subsection{(Further references not yet put in.)}
%
%Increases in car ownership:
%\begin{itemize}
%  \item Statistics on licensed vehicles in Great Britain, 1994 to 2013:  (Table VEH0102, graph included with the excel file) \url{https://www.gov.uk/government/statistical-data-sets/veh01-vehicles-registered-for-the-first-time}.  In 1994, 21 199 200 cars were licensed.  In 2013, 29 140 900 cars were licensed, an increase in amount by about 37\% of the 1994 figure.
%\end{itemize}
%
%\bigskip
%
%\noindent Pressure on parking facilities:
%\begin{itemize}
%  \item Newspaper article on the shortage of parking spaces in some Australian suburbs:  \url{http://www.dailytelegraph.com.au/newslocal/city-east/pressure-on-car-parking-spaces-in-potts-point-darlinghurst-and-elizabeth-bay-impacted-by-car-sharing-and-parking-permit-restrictions/story-fngr8h22-1227029587187?nk=be020414d91f4bfad788ba5547f430b8}
%\end{itemize}
%
%\bigskip
%
%\noindent Schemes:
%\begin{itemize}
%  \item encouraging the use of public transport, and park and ride, e.g. \url{http://www.luas.ie/park-and-ride/}
%  \item car sharing, e.g. \url{http://www.estates.salford.ac.uk/cms/resources/uploads/files/Car%20Share%20Scheme%20Frequently%20Asked%20Questions.pdf
%      }
%  \item renting parking space, and parking space sharing web sites, e.g. (1) \url{https://www.parkme.ie/}; (2) \url{https://www.yourparkingspace.co.uk/}; (3) \url{http://www.thisismoney.co.uk/money/bills/article-1614219/Should-you-rent-out-your-parking-space.html}
%  \item congestion charges
%  \item resident parking permits, e.g. \url{http://www.coventry.gov.uk/info/117/parking/708/residents_parking_schemes} (this one still doesn't guarantee that you'll get a parking space, you simply are able to compete)
%\end{itemize}

\section{Prior Work}
\label{sec:priorwork}

Within the research community, questions concerning how to manage parking space supply-demand mismatch are actively being investigated from a variety of different angles.  One important aspect concerns the delivery of up-to-date, accurate, real-time information to parking systems in order to achieve greatest system efficiency.  A car park monitoring and vacancy detection system is discussed in \cite{blumer2012} and consists of a camera and image analysis system that reports statistics.  The aim for such systems is to provide timely information on available parking spaces and minimize search time.
ParkYa is an early-stage company offering a service that signposts parking locations and allows users to pay for their parking through a smartphone application.  Driving directions to the parking spot are additionally available. In the near future, such services should benefit greatly from incorporating increased amounts of real-time data, concerning up-to-date details on current parking vacancies, and predicted travel times given current road conditions.  The system of \cite{rico} monitors and reports the state of availability of parking spaces and additionally uses context information generated by the city and its citizens to provide accurate responses to demand.\newline

A predictive approach is proposed in \cite{caliskan} and further studied in \cite{klappenecker}. The authors develop a method in which car parks are able to count the number of arriving and departing cars, and report these to participating cars. This allows them to predict the likelihood of a parking space being available at the estimated time that the car will arrive there. This work uses ideas from queueing theory to predict the occupancy upon arrival, with car parks being modelled as single server queues with a Poisson arrival process and exponentially distributed service times. It should be mentioned that this significant reduction in requirements by using a stochastic approach comes at the cost of certainty for the customers. The lack of a reservation system makes it possible that customers arrive to a fully occupied car park. The main drawback of the approach is, however, that ultimately the customers will want to use the information to make a decision whether to try their luck and drive to the car park or to go somewhere else. Accordingly, there is feedback embedded in the system which needs to be taken into consideration; namely, when users choose to drive to a car park based on the predictions made, they then affect the arrival process, rendering the model and predictions no longer valid. This latter issue is addressed in \cite{schlote}, where multiple car park load balancing is achieved using randomisation based on congestion signals from the participating car parks. Significant extensions of this work are given in \cite{ijc15}.

\bigskip

Certainly, access to real-time information, on its own, is not enough.  Given the delay between an user identifying and driving to a parking space, and the fact that other users may also compete for the same parking space, a successful parking system would benefit from either intelligent assignment or reservation capability to reduce localized congestion \cite{geng}.  The reservation system is championed in \cite{barone}, where an architecture for an intelligent parking assistant is proposed as a public car park management solution.
In \cite{geng}, the parking problem is viewed as a dynamic resource-allocation problem. Similarities to problems in communication networks are drawn, for which a host of tools and methods have been developed over the last decades. An online reservation system is proposed, where cars communicate their parking requirements and are assigned a parking space, which is then reserved and cannot be used by any other vehicle. A similar approach is proposed in \cite{teodorovic}, albeit with a different assignment routine, and also allows the user to specify the price that he is willing to pay. The main focus of the paper is revenue maximisation, but it is also possible to achieve other goals, such as reducing traffic levels or ensuring some sort of fairness between users from different social classes. 
% It is also concluded that the optimal assignment strategy depends on the vehicle arrival process.

\bigskip

In \cite{dellorco2003}, an agent-based simulation model for parking facilities management was developed, the goal of the tool being to provide better understanding of how complex urban traffic systems evolve under different parking pricing strategies and different levels of parking enforcement.  Agents included drivers, parking authorities, law enforcement and city government.  Models for cruising-for-parking behaviour were developed in \cite{kaplan}.
In \cite{king}, the Quality of Service (QoS) metrics is used to dimension a shared vehicle fleet such that satisfactory levels of access to the fleet were ensured to its pool of users.  While the topic of \cite{king} was not parking per se, the paper inspires some of the notions of our current work; that is, we consider the dimensioning of a special subset of parking spaces in our mini-city that we will call the {\em reserve}.  We discuss this problem in the following sections of the paper.

\section{Setting}
\label{sec:mini-city-parking}

\begin{figure}[h]
\centering
\includegraphics[scale=0.8]{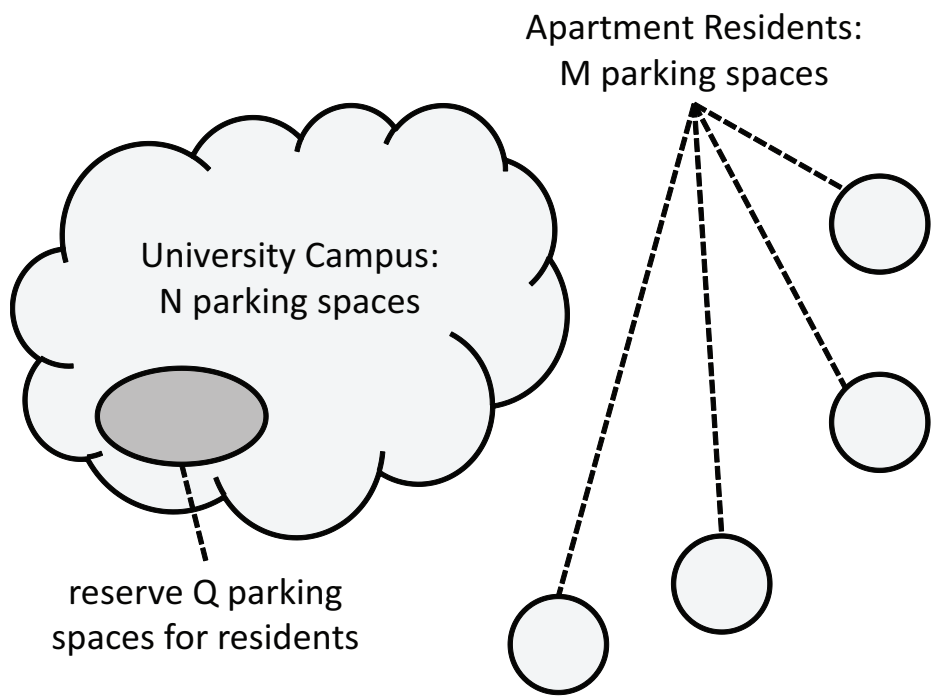}
\caption{The parking scenario: Premium spaces are those on the university campus, whereas secondary spaces are those belonging to the apartment buildings.}\label{problem_setup}
\end{figure}

We will consider the following scenario.  Imagine a university campus with a total of $N$ parking spaces, surrounded by a total of $M$ private parking spaces (e.g., from residential complexes), as illustrated in Figure \ref{problem_setup}.
A typical working day sees the university parking fill to capacity with vehicles belonging to students and staff, to the extent that the $N$ spaces cannot meet the demand.  Some campus arrivals thus have to search elsewhere.
At the same time, many nearby residents drive to work during the day and vacate their $M$ private parking spaces.
What we have is a wasted resource in one area (i.e., residential parking spaces) and a stressed resource nearby (i.e., the on-campus parking).  
We assume that a contract exists between the university and the \emph{landlords} or owners of these $M$ private parking spaces.
The contract stipulates that these landlords lease their driveways to the university during daytime
In the sequel, we will denote the time period in which the extra spaces are leased by $[0,W]$.
This can be possible because the landlords drive to work during this time period.

\bigskip

In our nomenclature, the premium resources are the on-campus parking spaces, while the secondary resources are the contracted residential parking spaces.
We will now consider how to design a parking system, bearing in mind the needs of the mini-city community, and the needs of the landlords.  We shall consider the following facets of the design problem:
\begin{enumerate}
  \item How do we accommodate the needs of parking space owners (landlords) given that situations may arise when landlords will need to return home during the contracted interval?  An example is scheduling a repairman to visit the home to fix a fitting. (See Section \ref{maths}.)
  \item How do we accommodate the needs of the landlords given that situations arise when misbehaving university members will not vacate the parking space outside of the contracted interval? (See Section \ref{maths}.)
  \item How do we allocate efficiently access to the premium spaces to the university community in a manner that preserves the privacy of individuals in the community? (See Section \ref{fair}.)
\end{enumerate}
As we shall see, we answer these questions in a stochastic framework by developing suitable Quality of Service (QoS) metrics to dimension the aforementioned parking system.  To address Items 1 and 2, we shall set aside a subset of $Q$ premium spaces as a reserve. The research question in this latter context then becomes how to use the QoS metrics to dimension the size of $Q$.  We shall address Item 3 by applying recent ideas from distributed resource allocation.

\bigskip

\textbf{Comment: (Generalization)}  The parking scenario presented in this paper is very specific, but our method can also be applied in the following general setting.  There are two resources, $A$ and $B$. These resources are atomic (indivisible), with each atom being allocated to one user. Each user has a preference over the two resources. Using resources $A$ and $B$ decreases a certain QoS metric. In the following main sections, we answer two questions:
\begin{itemize}
\item Given a QoS requirement, how much of resource $A$ and resource $B$ should we use?
\item Given an efficiency metric derived from the user preferences, how do we allocate the resources optimally among users in a setting where the allocation is repeated over multiple iterations?
\end{itemize}

\section{Dimensioning and Statistics}\label{maths}

In this section, we consider the problem of dimensioning $Q\in\{0,1,\ldots,N\}$ premium spaces as \emph{reserve} in order to provide sufficient Quality of Service (QoS) guarantees to the landlords of the secondary spaces.

\bigskip

Let $1,\ldots,M$ index the landlords of the parking spaces.
For each such landlord $i=1,\ldots,M$, we define a non-negative valued random variable $T_i$, which denotes the time at which he or she returns home and needs to get the parking space back.  Under normal circumstances, $T_i$ is greater than $W$, but on rare occasions a landlord may choose to come home early.

\bigskip

For simplicity, we assume that each parking space $i=1,\ldots,M$ has exactly one daytime user each day.
For each parking space $i=1,\ldots,M$, we also define a
non-negative valued random variable $A_i$, which denotes the departure time of the daytime user of the space $i$.  Under normal circumstances, $A_i$ is less than $W$, but we assume a small number of miscreants so that not all spaces are always vacated on time.

\bigskip

For convenience, recall that $[0,W]$ denotes the nominal rental window for every parking space. In other words, the landlord of space $i$ agrees to park only outside the interval $[0,W]$, whereas the daytime user of space $i$ agrees to park during the interval $[0,W]$ only.

\begin{defn}[Home-early and Overstay]
We define the following two events for each secondary parking space $i=1,\ldots,M$:
\begin{align*}
  E_i \triangleq \{ T_i \in [0,W] \} \cap \{T_i< A_i\},\\
  O_i \triangleq \{W < T_i\} \cap \{T_i< A_i\} .
\end{align*}
Each of these events represents an outcome where a landlord would like to use the space $i$, but cannot do so.
The \emph{home-early} event $E_i$ is due to the landlord needing the space during the day.
The \emph{overstay} event $O_i$ is due to the fact that the daytime user overstayed.
\end{defn}

For easy of presentation, we assume the following properties concerning the above random variables.

\begin{as}[For Simplicity]\label{as:1}
The random variables $\{T_1,\ldots,T_M\}$ are independent and identically distributed. Likewise,  $\{A_1,\ldots,A_M\}$ are also i.i.d.
Moreover, all these random variables are mutually independent.
Finally, we assume that all the distributions have densities.
\end{as}

\subsection{The Dimensioning Formulae}

We begin by quantifying the probability of an event $O_i$ in terms of the probability of daytime users overstaying.  Recall that the distribution of the random variable $A_i$ characterises the probability of the daytime user overstaying in space $i$.
First, observe that
\begin{align*}
  \P(O_i) \leq \P(A_i > W).
\end{align*}
Next, we derive an exact expression for $\P(O_i)$ using the
independence assumption (Assumption~\ref{as:1}).

\begin{lem}[Probability of $O_1$]
Let $F_{A}$ denote
the probability distribution of $A_1$; and $F_{T}$ denote the
probability distribution of $T_1$.
Under Assumption~\ref{as:1}, we have
  \begin{align*}
   \P(O_1)
           &= \int_{a=W}^\infty (F_{T}(a) - F_{T}(W)) dF_{A}(a).
\end{align*}
\end{lem}

\begin{rem}
  We estimate $F_{A}$ from data in the following section.
\end{rem}

\begin{proof}
Observe that
\begin{align*}
   \P(O_1) %&= \P\left( \{W < T_1< A_1\} \cap \{ A_1 > W \} \right)\\
           &= \P\left( W < T_1< A_1 \right)\\
           &= \int_{a=W}^\infty \P\left( W < T_1< A_1 \mid A_1 =
             a\right) dF_{A}(a)\\
           &= \int_{a=W}^\infty \P\left( W < T_1< a \right) dF_{A}(a).
\end{align*}
The claim follows by definition of the distribution $F_T$ of $T_1$.
\end{proof}

Now we derive a formula that we can use to dimension the reserve parking space $Q$ from the contingent of $N$ premium spaces.  Recall that there are $Q$ reserve parking spaces (the reserve ``buffer'') set aside by the university.
In this section, we consider the probability $p(M,Q)$ of the event that more than $Q$ spaces are needed to accommodate landlords needing the reserve buffer during daytime $[0,W]$:
\begin{align*}
  p(M,Q) = \P\left(\sum_{i=1}^M 1_{E_i} > Q \right),
\end{align*}
where $1_{E_i}$ denotes a Bernoulli random variable taking value 1 when event $E_i$ occurs, and value 0 otherwise.

\bigskip

First, we characterise the probability of the event $E_1$ in terms of
the probability distributions of $T_1$ and $A_1$.

\begin{lem}[Probability of $E_1$]
Let $F_{T}$ denote
the probability distribution of $T_1$.
Let $F_{A}$ denote
the probability distribution of
$A_1$. Under Assumption~\ref{as:1}, we have
\begin{align*}
  \P\left( E_1 \right) = \int_{t=0}^W (F_{A}(W) - F_{A}(t)) dF_{T}(t)
+ F_{T}(W) (1-F_{A}(W)).
\end{align*}
\end{lem}

\begin{rem}
$F_{T}$ and $F_{A}$ are estimated from data in
the next section.
\end{rem}

\begin{proof}
Let $\phi \triangleq \P\left( E_1 \right)$.
  Observe that
\begin{align*}
  \phi &= \P\left( \{T_1< A_1\} \cap \{ T_1 \in [0,W] \} \right)\\
       \mbox{(by non-negativity of $T_1$)}\quad&= \P\left( \{T_1 < A_1\} \cap \{ T_1 < W \} \right)\\
       &= \P\left( \{T_1 < A_1\} \cap \{ T_1 < W \} \mid A_1 \leq
         W\right) \P(A_1 \leq W)\\
&\quad+ \P\left( \{T_1 < A_1\} \cap \{ T_1 < W \} \mid A_1 >
         W\right) \P(A_1 > W)\\
       \mbox{(by Bayes' Rule)}\quad&= \P\left( T_1 < A_1 \mid A_1 \leq
         W\right) \P(A_1 \leq W)
+ \P\left( T_1 < W \right) \P(A_1 > W)\\
       &= \P\left( T_1 < A_1 \leq
         W\right)
+ \P\left( T_1 < W \right) \P(A_1 > W)\\
  &= \int_{t=0}^W \P\left( t < A_1 \leq
         W\right) dF_{T}(t)
+ \P\left( T_1 < W \right) \P(A_1 > W)\\
  &= \int_{t=0}^W (F_{A}(W) - F_{A}(t)) dF_{T}(t)
+ F_{T}(W) (1-F_{A}(W)),
\end{align*}
which is the claim.
\end{proof}

As a corollary, we have the following expression for the probability
that setting $Q$ reserve spaces at the university is not enough.

\begin{cor}[Probability that $Q$ reserve spaces are not enough]
  Let $\phi \triangleq \P\left( E_1 \right)$.
Under Assumption~\ref{as:1}, $p(M,Q)$ is a random variable
entirely characterized by $\phi$:
\begin{align}
  p(M,Q) = \sum_{k = Q}^M \binom{M}{k} \phi^k (1-\phi)^{M-k}.\label{eq:p}
\end{align}\label{cor_5}
\end{cor}

\begin{rem}
Observe that this probability can be mitigated by increasing the parameter $Q$.
Given a QoS target $p(M,Q)$ Equation \eqref{eq:p} can be used to determine the corresponding value of $Q$ needed to achieve it. This is illustrated through an example in the next subsection.
\end{rem}

\subsection{Parking Data and Example}

In the previous section we have considered arbitrary probability distributions $F_T$, $F_A$ in the formulae derived.
In this section, we give estimates $\hat F_T, \hat F_A$ of these distributions using publicly available data.
Given samples $Z_1,Z_2,\ldots,Z_n$ from the distribution $F_T$,
the corresponding empirical distribution-estimate takes the form
\begin{align*}
  \hat F_T(z) = \frac{1}{n} \sum_{i=1}^n 1_{[Z_i \leq z]}.
\end{align*}

\begin{figure}[h]
\centering
\includegraphics[scale=0.5]{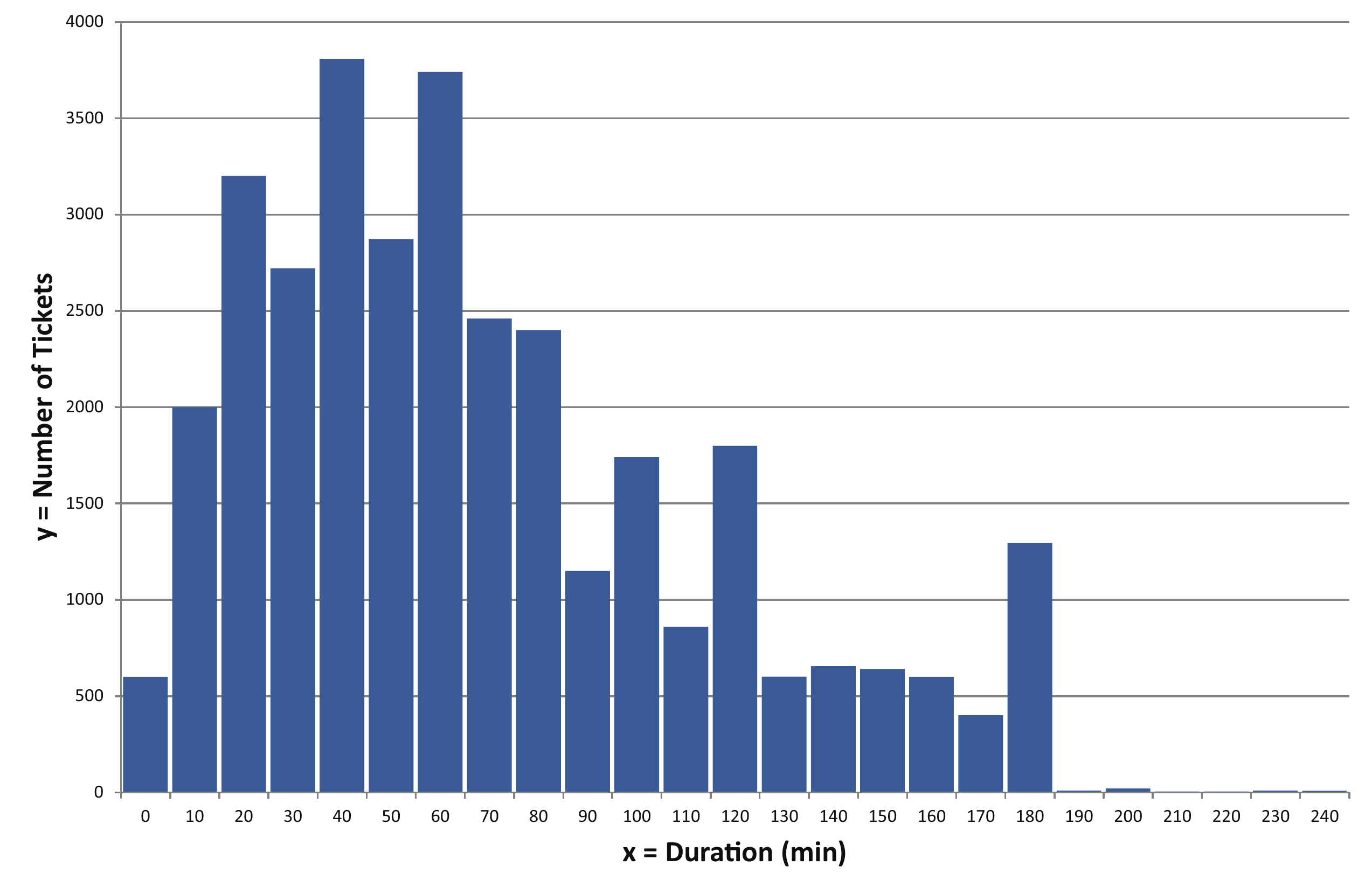}
\caption{Parking data.  (Source: Dublin City Council, 21 December 2013.)}\label{data}
\end{figure}

First, we estimate the distribution $F_A$.  Recall that $A_i$ is the random variable denoting the duration of use of the $i$th secondary parking space.  To derive an estimate, we use data on parking space utilisation collected in the city of Dublin.  Each data point corresponds to the time duration of one parking event. The histogram distribution of these durations is shown in Figure \ref{data}. Of course, we can obtain a better estimate of the distribution of parking usage in a university campus if we have access to more particular data.
Based on Figure~\ref{data}, in order to simulate the fact that 5\% of users of secondary parking spaces overstay, we set $W=170$ for our example.

% \begin{table}[h!]
% \caption{Probability distribution function of $A_i$.  For its construction, we map the qualitative properties of Figure \ref{data}.  We presume that the last approx. 5\% of the data mapped from Figure \ref{data} corresponds to parking overstay and hence take $W=170$ for our example.}
% \begin{center}
% \begin{tabular}{c|c|c|c|c|c|c|c|c|c|c|c|c|c|c|c|c|}\hline
% Time [min] &
% 0 & 
% 10 & 
% 20 & 
% 30 & 
% 40 & 
% 50 & 
% 60 & 
% 70 & 
% 80\\
% \hline
% Probability mass
% & 0.017
% & 0.059
% & 0.095
% & 0.080
% & 0.113
% & 0.085
% & 0.111
% & 0.073
% & 0.071\\
% Time &
% 90 & 
% 100 & 
% 110 & 
% 120 & 
% 130 & 
% 140 & 
% 150 & 
% 160 & 
% 170 \\
% Probability mass & 0.034
% & 0.051
%  & 0.025
%  & 0.053
%  & 0.017
%  & 0.019
%  & 0.019
%  & 0.017
%  & 0.011\\
% Time & 
% 180 &
% 190 & 
% 200 & 
% 210 & 
% 220 & 
% 230 & 
% 240 \\
% Probability mass
%  & 0.038
%  & 0.000
%  & 0.000
%  & 0.000
%  & 0.000
%  & 0.000
%  & 0.000
% \end{tabular}
% \end{center}\label{FAi}
% \end{table}

%\begin{figure}[h]
%\centering
%\includegraphics[scale=0.5]{FAi}
%\caption{Probability distribution function of $A_i$.  The function is constructed such that 5\% of the data mapped from Figure \ref{data} corresponds to parking overstay, and the area under the entire curve is equal to one.}\label{FAi}
%\end{figure}

% \bigskip

% \textbf{Comment:} In effect, we have artificially created a distribution for $A_i$. Of course, this could easily be measured by university authorities.

\bigskip

\begin{figure}[h]
\centering
\includegraphics[scale=0.5]{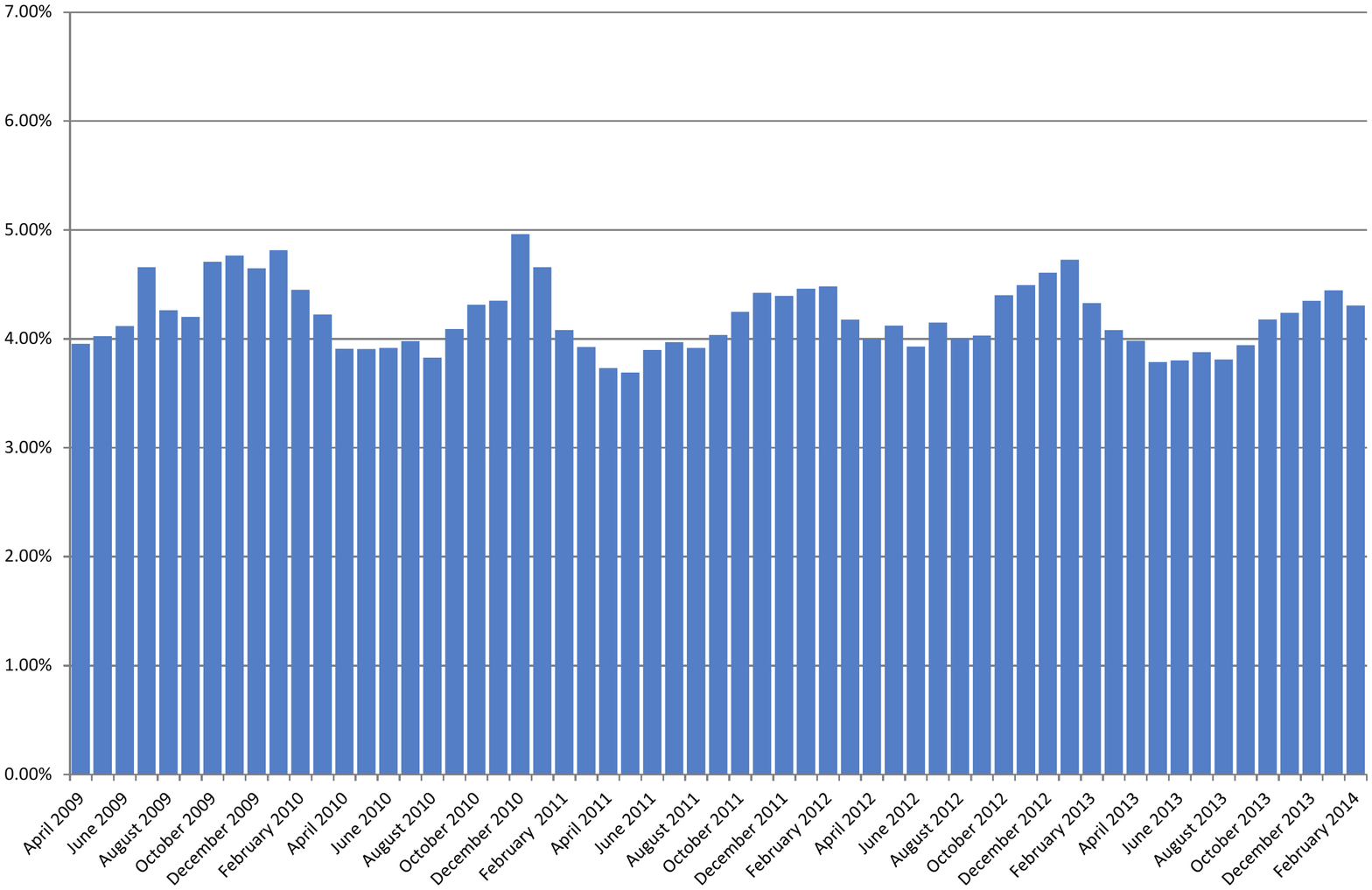}
\caption{Monthly sickness absence rates.  (Source: Health and Social Care Information Centre.)}\label{data2}
\end{figure}

Next, we estimate the distribution $F_T$ for $\{T_i\}$.  This distribution accounts for landlords who do not leave home and who arrive home normally after working hours.  There are many reasons for landlords to return home early, or not leave home at all.  For simplicity, we estimate the frequency of such days with the number of sick days of NHS staff in England over the period from April 2009 to February 2014.  This data is presented in Figure \ref{data2}.  For simplicity, we use the average sickness absence rate to estimate 
the probability of $T_i = 0$, and assume that $T_i = W$ otherwise (cf. Table~\ref{FTi}).

\begin{table}[h]
\caption{Probability distribution function of $T_i$.  For simplicity, we assume that the function is constructed so that 4.2\% of landlords remain at home due to sickness, while the other 95.8\% return home exactly at time instant $W$.}
\begin{center}
\begin{tabular}{c|c}\hline\hline
Random variable $T_i$ & Probability\\
\hline
$T_i =0$ & 0.042\\
$T_i =W$ & 0.958\\
\hline
\end{tabular}
\end{center}\label{FTi}
\end{table}

% go to doctor during day...uniform return home time

\bigskip

\begin{figure}[h]
\centering
\includegraphics[scale=0.6]{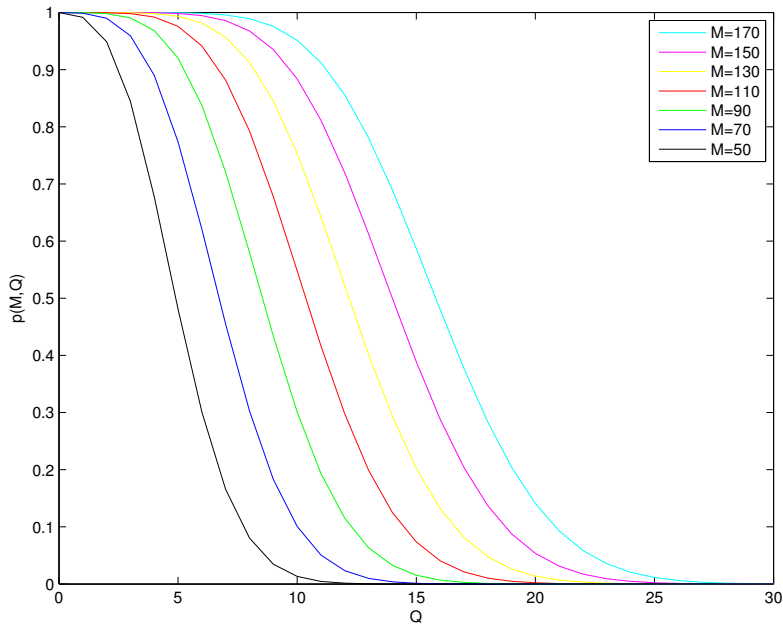}
\caption{Probability $p(M,Q)$ that $Q$ reserve spaces are insufficient when $M$ secondary spaces are contracted }\label{dimensioning_example}
\end{figure}

Finally, we use Corollary \ref{cor_5} to perform a dimensioning exercise based on our data.  Figure \ref{dimensioning_example} illustrates the probability $p(M,Q)$ that $Q$ reserve spaces are insufficient when $M$ secondary spaces are contracted. The probability $p(M,Q)$ eventually falls exponentially fast versus $Q$.  The dependence of $p(M,Q)$ on $M$ is more subdued. In other words, for a fixed value of $p(M,Q)$, a linear increase in the number of secondary spaces only requires a logarithmic increase in the number of reserve spaces.

\section{Efficient Allocation of Premium Spaces}
\label{fair}

In this section, we consider the point of view of the user (i.e. a staff member or student) who requires parking at the university.  Typically, such a user would purchase a monthly or yearly parking ticket.  This ticket then provides them with the opportunity to compete for a university parking space.  However, such a ticket does not necessarily guarantee them a parking space on the university grounds.  That is, if they arrive to the university ``too late'' on any given day, the car parks on-campus may already be full.  This is because more parking tickets are sold than there are actual spaces to park.

Such ``first come, first served'' systems can be inefficient.  Consider the example of a studying mother or father who, even though they paid the same amount for a monthly or yearly parking ticket as everyone else, must first always take their children to school in the morning.  They thus always arrive to the campus a little bit later than everyone else and have a greater chance of missing out on an on-campus parking space.  We are interested in a scheme that offers equality in regards to access to the premium, on-campus parking spaces for all university users over the long-term period of validity of their parking ticket.

For our scheme, we suppose that the university has been able to obtain enough apartment building parking spaces such that $(N-Q) + M$ is greater than or equal to the number of university users requiring a parking space on any given day.  The problem then is how to allocate the premium and secondary parking spaces to users efficiently over time.  How do we ensure that one campus user does not have to park off-campus for the majority of the time while another user is almost always offered a car space on the university grounds?  Our solution to this problem is motivated by additive-increase/multiplicative-decrease (AIMD) optimisation algorithms \cite{markovopt}.

%E.g. of fair scheme is buying stocks.  Company does well, everyone benefits.  Instead of paying fees to parking operator, drivers buy shares in parking and use dividends to pay for parking spot.  I.e. you buy a spot and generate revenue when people use it instead of you.  If high demand, the price goes up.

%Let the state $x_{i}$ denote the success of an individual campus goer at being

\bigskip

\subsection{Algorithm}\label{sec:fair-access-premium}

We now consider the problem of providing efficient access to the available premium spaces
at the mini-city. These spaces refer to the $N-Q$ spaces available after allocating $Q$ of the spaces as reserve according to the previous section. For simplicity of exposition, we shall assume that $Q = 0$ in this section, but the results generalize in a straightforward fashion.
% Recall our motivation is based on the observation that
% many services, not only parking, are based on users paying a fee to
% access, some limited resource. For example, many university campuses
% require students and staff to pay an annual parking fee. Maynooth
% University is an example of one such university \cite{} {\bf (Wynita:
%   please give reference to Maynooth packing policy)} - and there and many
% others. Similar services are also provided in a similar manner: seats on a
% train; access to an overhead bin in an aircraft; as well as access to
% shared services in clubs and in other societal groupings (car sharing is
% one such example).  In almost all examples of such services, despite
% appearance, payment entitles users to compete for the resource, and does
% not guarantee access. Often, users arriving first are almost always
% preferred. Our objective is now to attempt to address this problem by
% suggesting mechanisms to guarantee fair access. 
We introduce and recall the following notation to describe the problem data and the
variables used in our algorithm.
\begin{itemize}
\item[$N\\$ :] an integer denoting the number of available premium parking spaces.
\item[$n$ :] an integer denoting the number of users wishing to avail of premium parking spaces. We assume that $n > N$.
\item[$k$ :] denotes discrete time, $k=0,1,2,3,...$. In our interpretation
  this corresponds to the number of days the system is operating. For convenience, we assume spaces are assigned on a per-day basis, but the 
general principles of the algorithm do not depend on this assumption. 
\item[$X_i(k)$:] This is a state variable associated with the $i$th user. It takes the value $1$ if this user is given access to a premium parking on the $k$th day and zero otherwise.
\item[$\overline{X}_i(k)$:] This is a average access for the $i$th user up to the $k$th day:
\begin{equation*}
\overline{X}_i(k) = \frac{1}{k+1} \sum_{j=0}^{k} X_i(j) \,. 
\end{equation*}
\end{itemize}

It is possible to formulate the premium parking space allocation problem in several
ways. For example, one could require that that the long-term average
admission to the premium parking space is equal for all users, i.e. for
all $i,j = 1,\ldots,n$
\begin{equation*}
    \lim_{k\to \infty} \overline{X}_i (k) -\overline{X}_j (k) = 0 \,.
\end{equation*}
This assumes that all users are equal in the desire to access the premium
parking space. Here we follow a more general approach and assume that each
user $i$ has a cost function $f_i: [0,1] \to \mathbb R$. For a frequency $z\in [0,1]$ of premium space allocation, the cost $f_i(z)$ represents the monetary inconvenience cost that user $i$ experiences from $z$. This function specifies the
priority that this user is assigned. It could represent the amount a user
is willing to pay, or it could be related to the number of passengers
carried by this users, or the access that this user has to public
transportation (meaning that users with fewer possibilities for
alternative transport should have
prioritised access to parking spaces). Given these individual cost
functions, our aim is to design a system that achieves overall minimal cost
for the group. 
We formulate the optimal allocation of resources as a minimization problem:
\begin{eqnarray}
\label{eq:optprob}
\minimize_{z_1,\ldots,z_n \in \R} \,\, &\sum_{i=1}^n f_i (z_i)\\
\mbox{subject to}\quad& \sum_{i=1}^n z_i = N,\nonumber\\
& z_i\geq 0, \quad i=1,\ldots,n.\nonumber
\end{eqnarray}

Our proposed simple algorithm for solving the parking allocation problem can be summarized as follows. We assume that every event each user
is allocated access to the premium spaces by tossing a coin. For example, one embodiment of this idea is to use
a smart-phone application. 
More specifically, each user is assigned a cost function by a government authority.
For example, this could be based on vehicle class, disability, need for childcare, etc. 
On every subsequent day, the application outputs a variable that determines whether the user is assigned a premium space.
Access to the space is assigned to the $i$th user
with the following probability:
\begin{equation}
p_i(k) 
\triangleq \P(X_i(k) = 1) = \Gamma(k)\ \frac{\overline{X}_i(k)}{f_i'(\overline{X}_i(k))}\,.
\end{equation}
A central authority that owns
the parking spaces, calculates $\Gamma(k)$ based on past utilization and broadcasts this scalar to all participating vehicles. The following example
illustrates the performance of this simple algorithm.\newline

We will assume that the functions $f_i$ are continuously
differentiable and strictly
convex; so that in particular the optimal point $z^*\in \R^n$ satisfying
the constraints is unique. 
Furthermore, we introduce an assumption which ensures that the
optimal point $z^*$ has only positive entries. This assumption also
guarantees that the algorithm we will describe is well defined for every user.\newline

We wish to control the access to the premium space in such
a way that the average utilization for each user approaches the optimal
value $z^*$, i.e. for large $k$ we want to achieve 
\begin{equation*}
    \overline{X}_i (k) \approx z_i^* \,.
\end{equation*}
subject to the (loose) capacity constraint $\sum_{i=1}^n {x}_i \approx N$.
That is, all premium spaces are occupied, on average. Further, we wish to
do this in a manner that preserves the privacy of users. That is, we do
not wish to reveal $\overline{X}_i$ and $f_i$ to any other users in the course of the optimization. Finally, the necessary
communication between the users should be minimal so as not to create a
communication overhead that would be hard to sustain in an uncertain
environment where users cannot be expected to participate at all times.

 In what follows the mechanism
for preserving privacy is to develop a distributed algorithm.  We loosely follow the ideas in \cite{markovopt}, where a
distributed stochastic algorithm is presented which
guarantees that the average utilization variables $\overline{X}_i (k)$
converge to the optimal points $z_i^*$\footnote{As the algorithm is
  stochastic the convergence holds with probability $1$, which is also
  called almost sure convergence in a stochastic context.}.
The algorithm
presented here extends the ideas of \cite{markovopt}, however, because we
wish to ensure in addition that the instantaneous utilization variables
$x_i(k)$ sum to $N$, or at least to a value close to it.
Moreover, the resource to be allocated in our setting is atomic, as opposed to arbitrarily divisible.
These differences require
substantial changes to the algorithm presented in \cite{markovopt}.

At each time $k$, each user $i$ determines a probability $p_i(k)$ and
sets
\begin{equation}
\label{eq:alg1}
    X_i(k+1) = \left \{
      \begin{matrix}
          1 & \text{ with probability } p_i(k), \\ 
          0 &  \text{ with probability } 1-p_i(k),
      \end{matrix}
\right. \,
\end{equation}
where we note that all users make this probabilistic choice independently of
other users or previous decisions. 
The evolution of the probabilities is governed by the equation
\begin{equation}
    \label{eq:alg2}
    p_i(k) = \Gamma(k) \frac{\overline{X}_i(k)}{f_i'(\overline{X}_i(k))} 
\end{equation}
Note that each user $i$ can determine its own probability with the exclusive
knowledge of its own past utilization $\overline{X}_i(k)$ and cost function $f_i$. No information
from other users is required. The scalar $\Gamma(k)$ is a network wide
constant determined by the central agency and broadcast to all users.
Here, $\Gamma(k)$ is chosen such that $p_i(k) \in (0,1)$
for all $i=1,\ldots,n$ and all $k \in \N$. It is determined in a
time-varying manner as it also influences the demand for premium
spaces. Specifically, if at a certain time $k$ each $p_i(k)$ is fixed then
the expected utilization of the premium spaces is
\begin{equation}
    \label{eq:expectload}
    \mathbb{E} \left( \sum_{i=1}^n X_i(k+1) \right) = \sum_{i=1}^n p_i(k) =
    \Gamma(k) \sum_{i=1}^n
    \frac{\overline{X}_i(k)}{f_i'(\overline{X}_i(k))} \,.
\end{equation}
Moreover, the (random) instantaneous utilization $\sum_{i=1}^n X_i(k+1)$ is concentrated around the expected utilization by independence and Hoeffding's Inequality.
If we wish to ensure optimal utilization of the premium spaces and avoid
overbooking, the expected utilization should be below the number of
premium parking spaces. For instance, a standard deviation below this
number, depending on the desired quality-of-service metric. Denoting this
number by $N_E \leq N$, we will therefore adjust $\Gamma(k)$ so that the
expectation in \eqref{eq:expectload} tracks $N_E$. As the expectation is unknown, we use the
observed utilization as an estimator for this. Taking a simple error
regulation approach we thus arrive at 
\begin{equation}
    \label{eq:alg3}
    \Gamma(k+1) = \Gamma(k) + \alpha \left( N_E - \sum_{i=1}^n X_i(k) \right)\,.
\end{equation}

The overall system is now prescribed by the dynamics of $X_i$ given by
\eqref{eq:alg1}, the dynamics of $p_i$ given by \eqref{eq:alg2} and the
dynamics of $\Gamma$ as described in \eqref{eq:alg3}. 

A proof of the convergence of the algorithm is beyond the scope of
this paper. However, we present below an example to show its efficacy. 

\subsection{Example} 

\begin{figure}[h]
\centering
\includegraphics[height=2.5in]{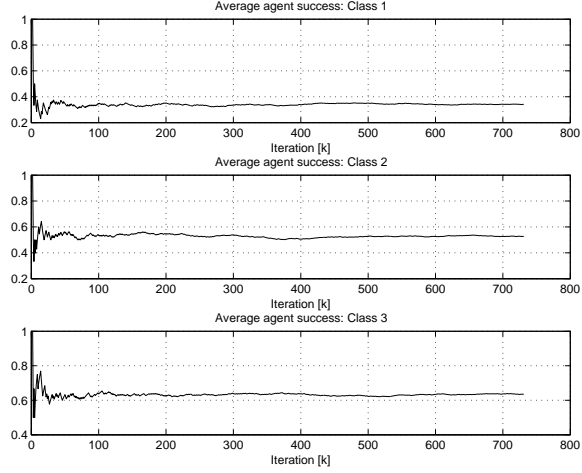}
\caption{Instantaneous allocation $\{X_i(k)\}$ for three users.}
\label{fig:ap}
\end{figure}

\begin{figure}[h]
\centering
\includegraphics[height=2.5in]{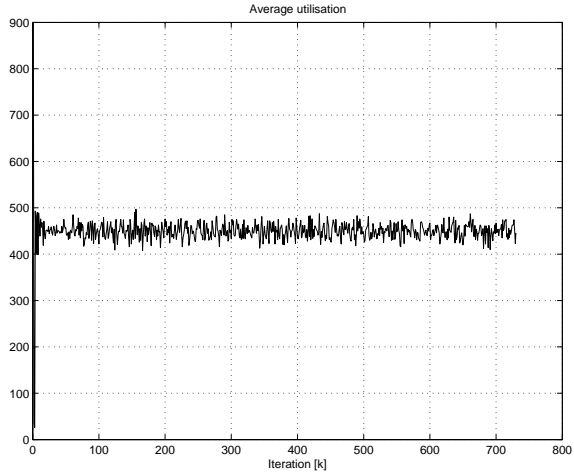}
\caption{Instantaneous premium space occupancy over time.}
\label{fig:pu}
\end{figure}

\begin{figure}[h]
\centering
\includegraphics[height=2.5in]{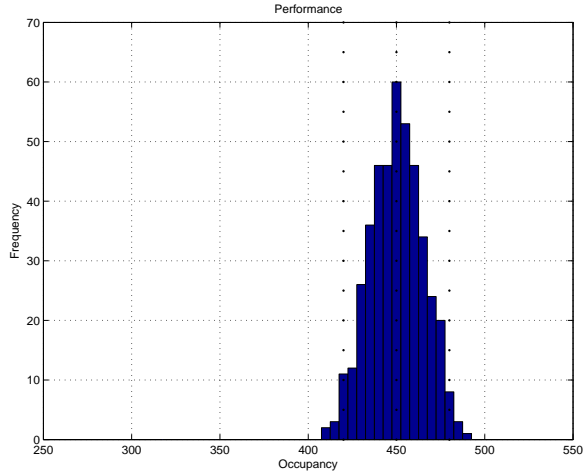}
\caption{Histogram of premium space occupancy over a horizon of 700 iterations.}
\label{fig:occ}
\end{figure}

\begin{figure}[h]
\centering
\includegraphics[height=2.5in]{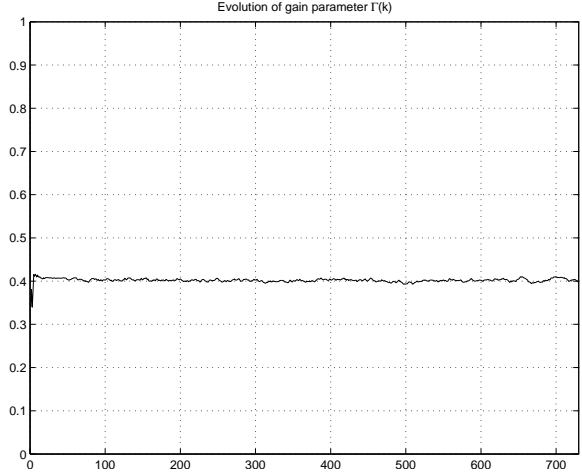}
\caption{Evolution of $\Gamma(k)$}
\label{fig:gain}
\end{figure}

\begin{figure}[h]
\centering
\includegraphics[height=2.5in]{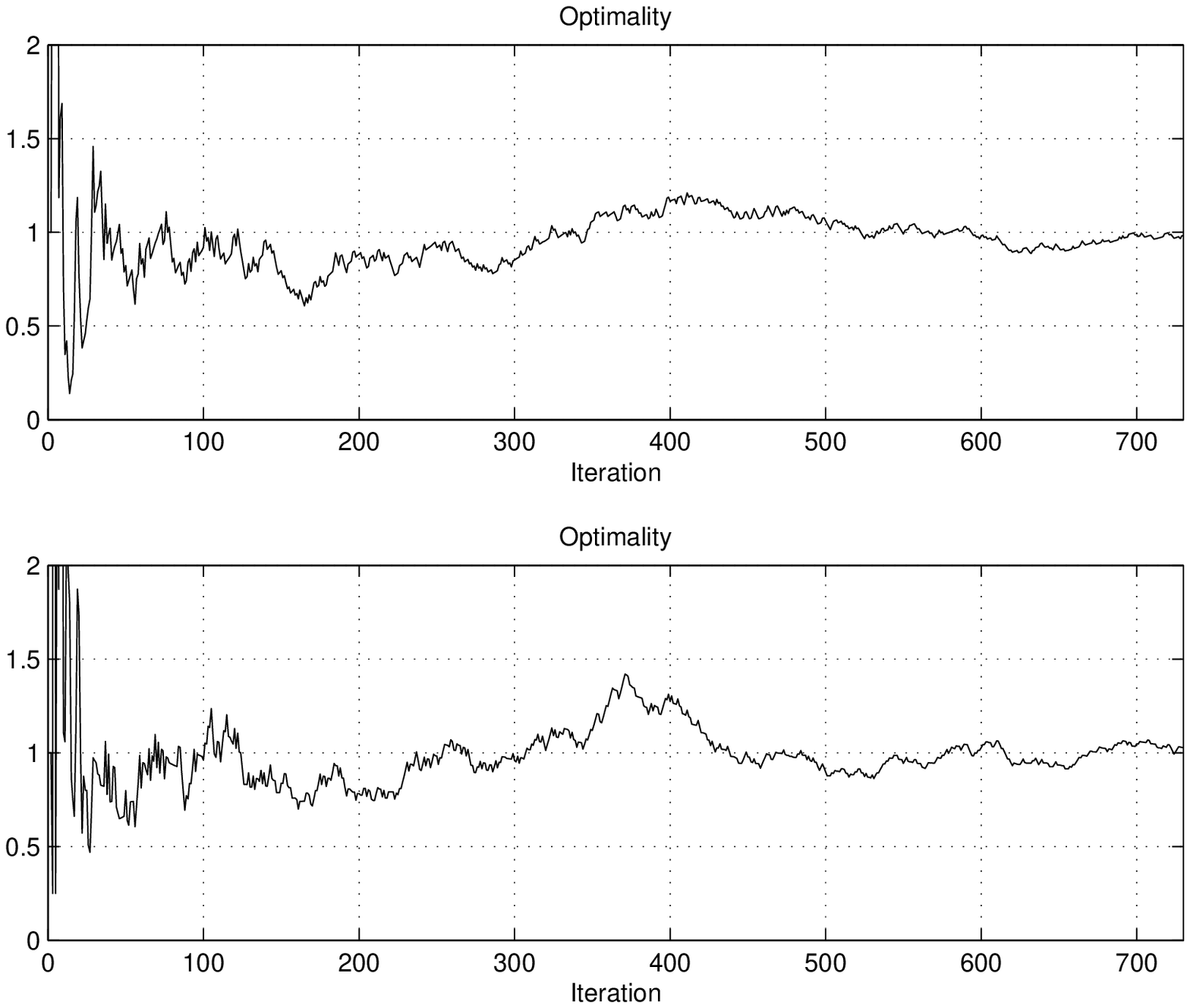}
\caption{Convergence to optimality. The marginal utilities of different users converge.}
\label{fig:con}
\end{figure}

%\end{example}

We simulate a population of 900 users competing for 450 premium parking spaces. Each evening, users
are assigned a parking space as described above (with the scalar $\Gamma(k)$ determined using a PI control). For simplicity, users have one of three cost functions: $f(z) = z^4/4$, $f(z) = z^6/6$, and $f(z) = z^8/8$, all of which are strictly convex.

%Figures \ref{fig:ap}--\ref{fig:occ} depict the performance of the algorithm. 
Figure \ref{fig:ap} shows the average utilisation
achieved for each class of vehicle. Figure~\ref{fig:pu} and~\ref{fig:occ} shows that the average utilisation of premium spaces is concentrated around the target utilisation of 450.
Figure \ref{fig:gain} shows the evolution of the controller gain and Figure \ref{fig:con} shows that the cost functions derivatives
do actually achieve consensus asymptotically.

\subsection{Intuition}

Here we will
briefly comment on the heuristic which suggests that the algorithm works
in a general setting. To this end we have to point out that a number of
assumptions are required.

First, in our derivation of the update rule \eqref{eq:alg2} we will rely
on an infinitesimal characterization of the optimal point. This requires
that the optimal point does not lie on the boundary of the constraint
set. Thus in the following we will assume that the optimal point $z^*$ satisfies
\begin{equation}
    \label{eq:ass1}
    z^*_i \in (0,1) \quad \forall \ i = 1, \ldots, n \,.
\end{equation}

Second, the gain parameter $\alpha$ in the feedback scheme \eqref{eq:alg3}
needs to be chosen in such a manner that the system is stable. Noting that 
by \eqref{eq:expectload} the expectation $\mathbb{E} \left( \sum_{i=1}^n
  X_i(k) \right)$ is a multiple of $\Gamma(k)$, we choose 
\begin{equation}
    \label{eq:2}
    \alpha \in \left( 0 , \left( \max_{x}  \sum_{i=1}^n
    \frac{x_i}{f_i'(x_i)}\right)^{-1}  \right)
\end{equation}
where the maximum with respect to $x$ is taken over the constraint
set. We then reformulate \eqref{eq:alg3} to obtain
\begin{equation*}
    \Gamma(k+1) = \left( 1 - \alpha \frac{\sum_{i=1}^n X_i(k)}{\Gamma(k)}\right) \Gamma(k) + \alpha N_E \,.
\end{equation*}
In the previous formulation, the bracket giving the factor of $\Gamma(k)$ needs to be in the interval $(0,1)$ sufficiently often to ensure stable behaviour of the system and thus the desired tracking property. The choice of $\alpha$ in \eqref{eq:2} is made to ensure that the desired stability property holds.
This
ensures convergence of the sequence $\Gamma(k)$.

To argue for the convergence of $\overline{X}$ we follow the heuristic
presented in \cite{markovopt}. The optimisation problem is a constrained
problem for which we have assumed that the boundary conditions $z_i \geq
0$ are not active. It is then classical to
introduce the Lagrange parameter $\mu \in \R$ and consider the Lagrangian
\begin{eqnarray}
    H(z_1,...,z_n, \mu) = \sum_{i=1}^n f_i  (z_i) - \mu \left(\sum_{i=1}^n z_i-N\right).
\end{eqnarray}
From the Karush-Kuhn-Tucker (KKT) conditions \cite[Section 5.5.3]{boyd2004convex}, the
following necessary and sufficient condition for optimality can be
obtained by setting all partial derivatives to zero. In the optimal point $z^* \in \R^n$, $\mu^*  \in \R$ we have
\begin{eqnarray}
\label{eq:KKT}
\mu^* = \frac{\partial f_i }{\partial z_i}(z_i^*) \;\; \quad \forall \; i= 1,\ldots,n.
\end{eqnarray}
In other words, the system is at optimality when the derivatives of the
cost functions are in consensus. This observation is at the core of the
choice of the probability functions in \eqref{eq:alg2} as we now explain.

Assuming that the system described by \eqref{eq:alg1}, \eqref{eq:alg2},
and \eqref{eq:alg3} has an ergodicity property we obtain the long-term average
for the access to the premium car park satisfies (in the limit)
\begin{eqnarray}
\overline{X}_i = \hat{p}_i
\end{eqnarray}
where $\hat{p}_i$ is the steady state probability that user $i$ has access to the resource.
Given this assumption, we claim that our choice of
 probability functions $p_i(k)$ is such that  the equation for the
steady state behaviour of our system is equivalent to the KKT
condition \eqref{eq:KKT}. By \eqref{eq:alg2} each user has a
probability of access that satisfies
\begin{equation}\label{eq:prob-set}
p_i(k) \propto  \ \frac{\overline{X}_i(k)}{f_i'(\overline{X}_i(k))}\,,
\end{equation}
where the constant of proportionality is the same for all users. Suppose
now that after a transient period we have $\overline{X}_i(k) \approx
\hat{X}_i$, the steady state long-run average which exists under our
ergodicity assumption. Then we can write $p_i(\overline{X}_i(k)) \approx \hat{p}_i$. Thus, in steady state,
\begin{equation}
p_i(k) \ \approx \frac{\Gamma}{f_i'(\overline{X}_i(k))} p_i(k)
\end{equation}
and so $f_i'(\overline{X}_i(k)) \approx \Gamma$ for all $i,j$. So in the
limit we expect $f_i'(\overline{X}_i(k)) \to \Gamma$ for all $i$. This
means that in the limit we converge to a situation where the derivatives
of the cost functions are in consensus. As we have seen from the KKT
conditions the consensus conditions characterizes the optimal point of the
optimization problem so that we obtain
\begin{equation*}
    \overline{X}_i(k) \to z_i^* \, \quad i=1,\ldots,n\,.
\end{equation*}
Thus the average utilization converges to the optimal point.

\subsection{Comments}

Several comments are merited.
\begin{itemize}
\item[(i)] This is a private algorithm. Users do not communicate with other users, nor do they reveal any state information or cost
information to the central authority.

\item[(ii)] This is a randomized algorithm. Thus we only have $\sum_i X_i(k) \approx N$. We argue that small fluctuations in the occupancy level
about the capacity is acceptable since some users may not use the premium spaces, even if allocated one. Furthermore,
in the event of all users availing of their allocated premium spaces, extra space is always available in the secondary spaces, and in the reserved spacing.
That is, the instantaneous allocation $\sum_i X_i(k)$ may exceed the capacity $N$ for some time steps $k$, in which case, some users will have to be re-routed to the secondary parking spaces.

\item[(iv)] The algorithm can be modified such that the parameter $\Gamma$ does not update at every iteration, but once in a while (i.e., from day to day).

\item[(v)] The proposed allocation scheme can be applied to more general problems of the form of \eqref{eq:optprob}. For example, the same problem arises also in assigning space in overhead bins on passenger planes, seats in trains, etc. In all of these examples, users pay a fee to compete for access to a limited resource. However, access to the resource is not guaranteed and in effect the resource is allocated in a first-come, first-served manner. From a societal viewpoint, customers with constraints, e.g. dropping children at school, are almost always disadvantaged by such schemes.

\end{itemize}

\section{Conclusion and Future Work}
\label{sec:concl-future-work}

In this work, we have modeled two aspects of parking system design and proposed two corresponding solutions.
First, we propose to complement primary or premium parking spaces with a large number of secondary parking spaces. These spaces are obtained through contracts for shared use that provide probabilistic QoS guarantees at the expense of a relatively small number of reserve premium spaces.
Secondly, we propose a distributed algorithm for repeated allocation of premium spaces, which guarantees that, over time, the users who benefit most from such spaces obtain them proportionally more frequently. The main feature of this algorithm is to keep private all information related to the benefit derived by individual users.

The next step in our research agenda is to deploy the described system on a university campus.
From an algorithmic perspective, it remains an open problem to 
quantify the rate of convergence of $X_i(k)$ and to situations where we have dynamic arrivals and departures throughout each day.
Another open problem is to quantify the probability of overshooting the capacity: i.e., the event $\sum_i X_i(k) > N$.
It is also useful to consider a version of the allocation algorithm for the case where a different subset of users competes in each day.

\end{document}